\documentclass{amsart}

\usepackage{epsfig}

\newcommand{\tr}{{\rm tr}}
\newcommand{\cL}{{\mathcal L}}

\newcommand{\cG}{{\mathcal G}}

\newcommand{\cI}{{\mathcal I}}

\newcommand{\bA}{{\hspace{-0.3pt}\mathbb A}\hspace{0.3pt}}
\newcommand{\bC}{{\hspace{-0.3pt}\mathbb C}\hspace{0.3pt}}
\newcommand{\bF}{{\hspace{-0.3pt}\mathbb F}\hspace{0.3pt}}
\newcommand{\bR}{{\hspace{-0.3pt}\mathbb R}\hspace{0.3pt}}
\newcommand{\bZ}{{\hspace{-0.3pt}\mathbb Z}\hspace{0.3pt}}

\newcommand{\g}{SL_2(\bC)} 
\newcommand{\Uh}{U_h(sl_2)}
\newcommand{\B}{\mbox{}_qSL_2}
\newcommand{\U}{U(sl_2)}

\newcommand{\ch}{\chi_\rho}

\newcommand{\lmk}{\mbox{}\hfill}
\newcommand{\rmk}{\hfill\mbox{}}

\newcommand{\lcr}{\raisebox{-5pt}{\mbox{}\hspace{1pt}
                  \epsfig{file=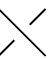}\hspace{1pt}\mbox{}}}
\newcommand{\ift}{\raisebox{-5pt}{\mbox{}\hspace{1pt}
                  \epsfig{file=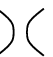}\hspace{1pt}\mbox{}}}
\newcommand{\zer}{\raisebox{-5pt}{\mbox{}\hspace{1pt}
                  \epsfig{file=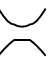}\hspace{1pt}\mbox{}}}

\newtheorem{theorem}{Theorem}

\title{Skein Quantization and Lattice Gauge Field Theory}
\makeatletter
\author{doug bullock}

\address{Department of Mathematics, Boise State University, Boise, ID
  83725, USA\\ email:  bullock@math.idbsu.edu}

\author{CHARLES FROHMAN}

\address{Department of Mathematics, University of Iowa, Iowa City, IA
  52245, USA\\ email: frohman@math.uiowa.edu}

\author{JOANNA  KANIA-BARTOSZY\'{N}SKA}

\address{Department of Mathematics, Boise State University, Boise, ID
  83725, USA\\ email: kania@math.idbsu.edu}

\thanks{The first author is partially supported by an Idaho SBOE
Specific Research Grant; the second and third by NSF-DMS-9204489 and
NSF-DMS-9626818.}

\makeatother

\begin{document}

\maketitle

\section{introduction}

The Kauffman bracket skein module is a deceptively simple
construction, occurring naturally in several fields of mathematics and
physics.  This paper is a survey of the various ways in which it is a
quantization of a classical object.

Przytycki \cite{P1} and Turaev \cite{T1} introduced skein modules.
Shortly thereafter, Turaev \cite{T2} discovered that they formed
quantizations of loop algebras; further work in this direction was
done by Hoste and Przytycki \cite{HPquant}.  We will look at some of
the heuristic reasons for treating skein modules as deformations, and
then realize the Kauffman bracket module as a precise quantization in
two different ways.

Traditionally, this is done by locating a non-commutative algebra that
deforms a commutative algebra in a manner coherent with a Poisson
structure.  The importance of the Kauffman bracket skein module began
to emerge from its relationship with $\g$ invariant theory.  It is
well known that the $\g$-characters of a surface group form a Poisson
algebra \cite{BG,Go}.  The skein module is the appropriate deformation.

The idea of a lattice gauge field theory quantization of surface
group characters is due to Fock and Rosly \cite{FR}. It was developed
by Alekseev, Grosse and Schomerus \cite{AGS} and by Buffenoir and
Roche \cite {BR1,Bu}.  We tie the approaches together by showing that
the skein module coincides with the lattice quantization.

\section{The Kauffman Bracket Skein Module}

\begin{figure}[t]
\centering
\epsfig{file=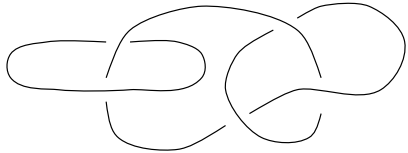}
\caption{Diagram of a two-component link.}
\label{link-diagram}
\end{figure}

Quantum topology began with the discovery of several new link
polynomials, the first and most well known being the Jones polynomial
\cite{jones1}, \cite{jones2}.  Many subsequent invariants arose from
alternative proofs of its existence.  The state sum approach
\cite{kauffman1} yielded an invariant known as the Kauffman bracket,
on which we will focus.  The Kauffman bracket is a function on the set
of framed links in $\bR^3$.  Since we will take a combinatorial view
throughout, one may as well think of a link as represented by a
diagram in $\bR^2$ (see Figure \ref{link-diagram}).

\begin{figure}[b]
\centering
\epsfig{file=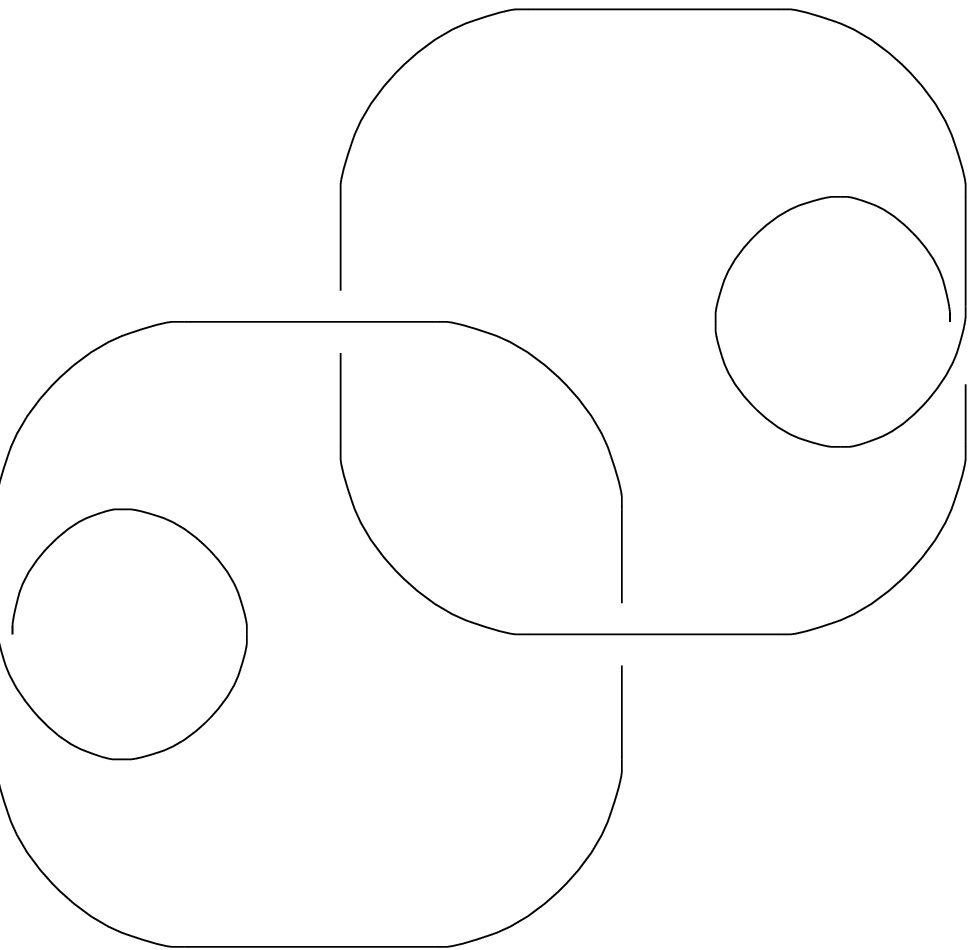,width=1in,angle=-45}
represents
\epsfig{file=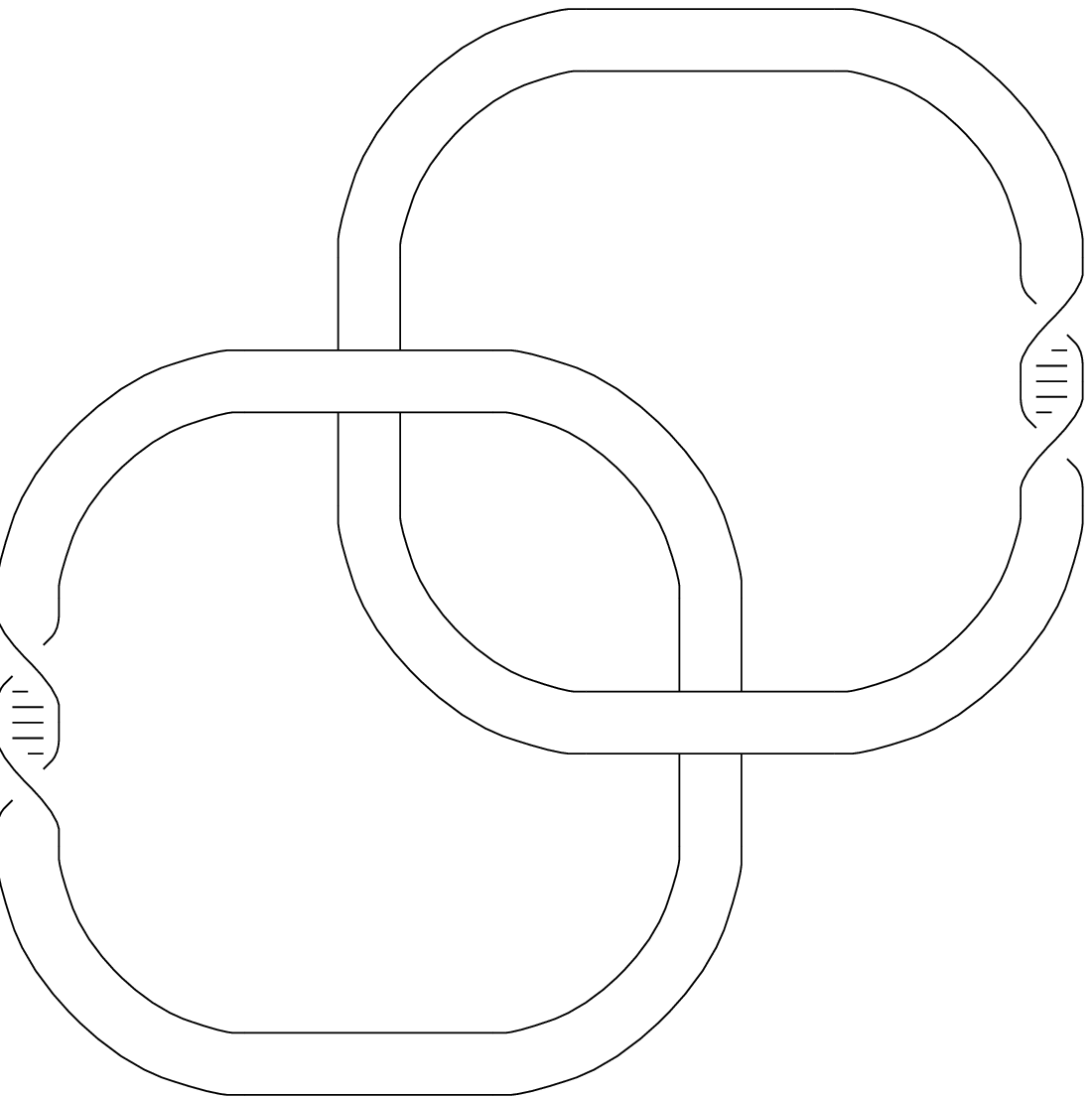,width=1in,angle=-45}
\caption{Diagram for a framed Hopf link}
\label{framed-hopf}
\end{figure}

A link is an embedding of circles, of which a diagram is a
particularly convenient picture.  Two diagrams represent the same link
if one can be deformed into the other.\footnote{These deformations may include some Reidemeister moves.}
In a framed link each circle
is actually the centerline of an embedded annulus.  Since we always
work with diagrams, it makes sense to assume the annulus lies flat in
$\bR^2$ as illustrated in Figure \ref{framed-hopf}.

The Kauffman bracket, $\langle\; \rangle$, takes values in the ring
$\bZ[A^{\pm1}]$ and is uniquely determined by the
rules:\footnote{Actually, Kauffman did not include the empty diagram
$\emptyset$, and his normalization was $\langle\text{unknot}\rangle =
1$.  Later, the quantum group approach \cite{KM,RT1} to knot
polynomials---and to some extent skein modules---indicated that
normalization at $\emptyset$ was preferable.}
\begin{enumerate}
\item $\langle\emptyset\rangle=1$,
\vspace{5pt}
\item $\displaystyle{\left\langle\lcr\right\rangle =
                   A \left\langle\zer\right\rangle +
              A^{-1} \left\langle\ift\right\rangle}$, and
\vspace{5pt}
\item $\displaystyle{\left\langle\bigcirc\right\rangle 
          = (-A^2-A^{-2})\langle\emptyset\rangle}$.
\end{enumerate}
The arguments of $\langle\; \rangle$ in (2) and (3) represent diagrams
which are identical except in a neighborhood where they differ as
shown in the formulas.  One evaluates the function on a diagram by
first applying (2) until no crossings remain, and then reducing each
diagram to a polynomial via (3) and (1).

\begin{theorem}[Kauffman]
The function $\langle\; \rangle$ is well defined. If $D_1$ and $D_2$
represent the same framed link, then $\langle D_1\rangle =\langle
D_2\rangle$.
\end{theorem}

Kauffman's construction is an example of a link invariant defined by
skein relations on the set of all diagrams.  Since skein relations are
defined only in small neighborhoods, the idea generalizes naturally to
spaces locally modeled on $\bR^3$.

The notion of a skein module of a $3$-manifold\footnote{A suggestion
appeared in Conway's treatment of the Alexander polynomial
\cite{conway}.} was introduced independently by Przytycki in \cite{P1}
and Turaev in \cite{T1}.  Roughly speaking, the construction consists
of dividing the linear space of all links by an appropriate set of
skein relations, usually the same as those known to define a
polynomial invariant in $\bR^3$.  We will give the explicit definition
for the Kauffman bracket skein module.

Let $\cL_M$ be the set of framed links (including $\emptyset$) in a
$3$-manifold $M$. Denote by $\bC\cL_M$ the vector space consisting of
all linear combinations of framed links.  Take $\bC\cL_M[[h]]$ to be
formal power series with coefficients in $\bC\cL_M$, and give it the
$h$-adic topology.\footnote{Skein modules were originally less
technical \cite{Lick2,P1,T1}.  Power series first appeared in
\cite{T2}.  Topological considerations were first addressed in
\cite{quant}.}  This is an example of a topological module (see
\cite{kassel} for a nice introduction), however, one may think of
$\bC\cL_M[[h]]$ as just the completion of a vector space with basis
$\cL_M$ and scalars $\bC[[h]]$.

Let $t$ denote the formal series $e^{h/4}$ in $\bC[[h]]$.  We define
the module of skein relations, $S(M)$, to be the the smallest subset
of $\bC\cL_M[[h]]$ that is closed under addition, multiplication by
scalars and the $h$-adic topology, and which contains all expressions
of the form
\begin{enumerate}
\item $\displaystyle{\lcr+t\zer+t^{-1}\ift}$, \quad and
\vspace{5pt}
\item $\bigcirc+t^2+t^{-2}$.
\end{enumerate}
As before, (1) and (2) indicate relations that hold among links which
can be isotoped in $M$ so that they are identical except in the
neighborhood shown.  The Kauffman bracket skein module is the quotient
\[K(M) = \bC\cL_M[[h]]/S(M).\]

This process can be mimicked for any choice of basis (oriented links,
links up to homotopy, etc.), any choice of scalars, any set of skein
relations, and with or without requiring topological completion.  The
resulting quotient is generically called a skein module.  For
instance, an older version of the Kauffman bracket skein module is
\[K_A(M)=\bZ[A^{\pm 1}]\cL_M / S(M),\]
with $t=-A$ in the skein relations, and without topology.  If
 $M=\bR^3$ (or $B^3$  or $S^3$)
the new version is just an outrageous way of expanding the Kauffman
bracket into a power series.

\begin{theorem}[Kauffman--Przytycki--B--F--K] $K(\bR^3) \cong \bC[[h]]$
via $L \mapsto \langle L\rangle_{A=-t}$.
\end{theorem}

On one level, $K(M)$ is a generalization of the Kauffman bracket
polynomial.  If $K(M)$ is topologically free (i.e. isomorphic to
$V[[h]]$ for some vector space $V$) then the isomorphism gives a power
series link invariant for each vector in a basis of $V$.  The
coefficients behave like finite type link invariants \cite{quant},
generalizing a well known property of the Jones polynomial expanded as
a power series \cite{birman-lin}.  In order to utilize the module in
this fashion, one would like to know that it is free, and what the
basis is; information that is decidedly difficult to come by.

The survey article \cite{HP3} contains a nearly complete list of those
manifolds for which the explicit computation has been done, the
exception being \cite{skein}.  These computations predate the
topological version of the module, but whenever $K_A(M)$ is free,
$K(M)$ is just its completion after substituting $A=-t$.

There is, however,  a  deeper  understanding of skein modules,
$K(M)$ in particular.  Przytycki often refers to skein theory as
``algebraic topology based on knots,'' alluding strongly to 
skein modules as a sort of non-commutative alternative to homology.
This is also reflected in the principle that loops up to homotopy
carry classical information, whereas knots up to isotopy carry quantum
information.  The notion that a skein module is a quantization or a
deformation of some kind can be made very explicit for $M=F \times I$,
$F$ being a compact oriented surface.

In this case, $K(F\times I)$  is an
algebra.  Multiplication of links in $\bC\cL_{F\times I}$ is by
stacking one atop the other; it extends obviously to $\bC\cL_{F\times
I}[[h]]$, and it is a simple matter to check that $S(F \times I)$ is
an ideal.  Crossings form a barrier to commutativity in
$\bC\cL_{F\times I}$, and, for most surfaces, the obstruction survives
in the quotient.\footnote{ The exceptions are planar surfaces with
$\chi(F) \geq -1$.}

It is possible for non-homeomorphic surfaces $F_1$ and $F_2$ to have
homeomorphic cylinders, $F_1\times I$ and $F_2\times I$. The
homeomorphism does not preserve the algebra
structure.\footnote{$F\times I$ in Theorem \ref{BP} is homeomorphic to
$\raisebox{-2pt}{\epsfig{file=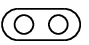,height=8pt}} \times I$, whose
skein algebra is commutative.} For this reason, it makes sense to
compress to notation into $K(F) = \bC\cL_F[[h]]/S(F)$.
 
For an example of $K(F)$ as a  ``deformation'', recall that the
commutative algebra 
of polynomials in three variables (over your favorite scalars) is
presented by
\[ \langle x,y,z\ |\ xy-yx=0, yz-zy=0, zx-xz=0 \rangle.\]  
 
\begin{theorem}[B--Przytycki]\label{BP} If $F$ is a  once punctured
torus, then $K_A(F)$ is presented by\footnote{The variables $x$, $y$
and $z$ are a meridian, a longitude and a slope one curve on $F$.}
\begin{align*}
\langle x,y,z\ |\ &Axy - A^{-1}yx=(A^2-A^{-2})z \\
                  &Ayz - A^{-1}zy=(A^2-A^{-2})x \\
                  &Azx - A^{-1}xz=(A^2-A^{-2})y \rangle.
\end{align*}
\end{theorem}

Theorem \ref{BP} can be thought of as a $1$-parameter family of
presentations,\footnote{Barrett \cite{barrett} showed that a spin
structure on $M$ induces an isomorphism $K_A(M)\cong K_{-A}(M)$.}
which reduces to the commutative polynomials if $A=\pm 1$.  Other
examples can be found in \cite{bptorus}.

Yet another way to see the module as a deformation is to let the
parameter $h$ go to $0$ (or $t$ to $1$, or $A$ to $-1$).  Formally,
this is achieved by passing to the quotient $K_0(F)=K(F)/hK(F)$. In
this case the skein relations would become
\begin{enumerate}
\item $\displaystyle{\lcr+\zer+\ift}$, \quad and
\vspace{5pt}
\item $\bigcirc+2$.
\end{enumerate}
Taken together these allow crossings to be changed at will and make
framing irrelevant.  Hence, the ``undeformed'' module is a commutative
algebra spanned by free homotopy classes of collections of loops.
The multiplication  in this algebra is commutative.

This is all quite heuristic for we lack a precise definition of
quantization or deformation. The next section will address this.  Even
so, if one can only understand the underlying commutative algebra as
an obvious quotient of the deformation, there is little here but
tautology.  We will close this section with an interpretation of
$K_0(F) = K(F)/hK(F)$ in terms of group characters.

Suppose that $G$ is a finitely presented group with generators
$\{a_i\}_{i=1}^m$ and relations $\{r_j\}_{j=1}^n$.  The space of $\g$
representations of $G$ is a closed affine algebraic
set.\footnote{Shafarevich \cite{shaf} is a good reference for the
algebraic geometry.}  You can view the representations as lying in
$\prod_{i=1}^m \g \subset \bC^{4m}$.  Each of the relations $r_j$
induces four equations from the coefficients of $r_j(A_1, \ldots
,A_m)=I$.  The zero set of these polynomials restricted to the variety
$\prod_{i=1}^m \g$ is the representation space.

We might naively try to construct the coordinate ring of the
representations as follows. Let $\cI$ be the ideal generated by the
equations $r_j(A_1, \ldots ,A_m)=I$ in the coordinate ring $
C[\prod_{i=1}^m\g]$.  Let $R(G)=  C[\prod_{i=1}^m\g]/\cI$.  The
problem is that $R(G)$ might have nilpotents (equivalently,
$\sqrt{\cI} \neq \cI$).  However, it was proved in \cite{LM} that
$R(G)$ is an isomorphism invariant of $G$. There is an action of $\g$
on $R(G)$ induced by conjugation in the factors of $\prod_{i=1}^m
\g$. The part of the ring fixed by this action is the affine
$\g$-characters of $G$, denoted $R(G)^{\g}$.  It, too, is an isomorphism
invariant of $G$ \cite{BH}. For our purposes, it suffices to define
the ring of $\g$-characters of $G$ to be\footnote{It is a deep result
of Culler and Shalen \cite{CS} that the set of characters of $\g$
representations of $G$ is a closed affine algebraic set.  Its
coordinate ring is $\Xi(G)$.}
\[ \Xi(G) = R(G)^{\g}/\sqrt{0}.\]
If $X$ is a manifold, we will write $\Xi(X)$ rather than
$\Xi(\pi_1(X))$.

The connection with 3-manifolds is quite simple (see
\cite{reps,estimate,isomorphism,PS1} for details).  Suppose that
$\rho:\pi_1(M) \rightarrow \g$ is a representation and $\ch$ is its
character.  Let $K$ be a loop, thought of as a conjugacy class in
$\pi_1(M)$.  Since the trace of a matrix in $\g$ is invariant under
inversion and conjugation, it makes sense to speak of $\ch(K)$
regardless of the choice of a starting point or orientation.  The loop
$K$ determines an element of $R(G)^{\g}$ by $K(\rho) = -\ch(K)$.  The
function extends to
\[ \Phi :K_0(M) \rightarrow R(G)^{\g}\]
by requiring it to be a map of algebras.  It is well defined because
the relations in $K_0(M)$ are sent to the fundamental $\g$ trace
identities:
\begin{enumerate}
\item $\tr(AB)+\tr(AB^{-1})=\tr(A)\tr(B)$, and 
\item $\tr(I)=2$. 
\end{enumerate}

It is shown in \cite{isomorphism} that the image of $\Phi$ is a
particular presentation of the affine characters \cite{GM}.  Sikora
\cite{sikora} has achieved this by directly identifying a version of
$K_0(M)$ with $R(\pi_1(M))^{\g}$ as defined by Brumfiel and Hilden.
Przytycki and Sikora \cite{PS1,PS2} have computed $K_0(M)$ for a large
number of manifolds, including $M=F\times I$, for which they can prove
it has no nilpotents.\footnote{They have also worked with various
scalars.  Certainly if the scalar ring has nilpotents then $K_0(M)$
does as well, and they have even located a nilpotent with scalar field
$Z_2$.  However, no nilpotents have ever been found with scalar ring
$\bC[[h]]$.}  Summarizing, we have a good idea of what $K_0(M)$ is in
general, and we know exactly what $K_0(F)$ is.

\begin{theorem}[B-Przytycki-Sikora]\label{BPS}
$\Phi :K_0(F) \rightarrow \Xi(F)$ is an isomorphism.
\end{theorem}

\section{Poisson Quantization of Surface Group Characters}

In the previous section we saw how a non-commutative algebra shrinks
to a commutative specialization for some particular value of a
deformation parameter. The formal definition of quantization reverses
this process. Beginning with a commutative algebra, one introduces a
parameter $h$, and a ``direction'' of deformation. The direction is a
Poisson bracket.

To make this precise, a commutative algebra $A$ is
called a Poisson algebra if it is equipped with a bilinear, antisymmetric map
$ \{\ ,\ \} : A \otimes A \rightarrow A$  which
satisfies the Jacobi identity:
\[ \{a,\{b,c\}\} +\{b,\{c,a\}\} + \{c,\{a,b\}\} = 0 , \]
and is a derivation:
\[ \{ab,c\} = a\{b,c\}+ b\{a,c\}, \]
for any $a,b,c \in A$.

A quantization of a complex Poisson algebra $A$ is a
$\bC[[h]]$-algebra, $A_h$, together with a $\bC$-algebra isomorphism,
$\Phi: A_h/hA_h \rightarrow A$, satisfying the following properties:
\begin{itemize}
\item   as a $\bC[[h]]$-module $A_h$ is topologically free
(i.e. $A_h\equiv V[[h]]$);
\item if 
$a,b \in A$ and $a',b'$ are any elements of $A_h$ with $\Phi(a')=a$
and $\Phi(b')=b$, then 
\[\Phi\left( \frac{a'b'-b'a'}{h}\right) = \{a,b\}. \] 
\end{itemize}

Hoste and Przytycki \cite{HPquant}, and Turaev \cite{T2} knew that
certain skein modules gave Poisson quantizations of various algebras
based on loops in a surface.\footnote{Their modules have a slightly
different flavor than the one defined here, both because topology is
not considered and because the scalars are not necessarily power
series.  Their definitions of Poisson quantization are analogously
distinct.  It is also interesting to note that their work predates the
appearance of quantum groups in low-dimensional topology.}  Since
$K(F)$ is topologically free (\cite{P1}, \cite{quant}), one can easily
see it as a Poisson quantization of $K_0(F)$ with the obvious bracket:
\[\{a,b\}=\text{lead coefficient of $a'b'-b'a'$ in $K(M)$.}\]
As noted in Section 3, however, understanding the Poisson algebra $K_0(F)$
as a formal quotient of $K(M)$ yields no new insight.
This is where character theory reenters.  

Since $\Xi(F)$ is the complexification of the $SU(2)$-characters of
$\pi_1(F)$, it has a Poisson structure given by complexifying the
standard one on $SU(2)$-characters \cite{Go,BG}.  Recall
(Theorem \ref{BPS}) that the algebra $\Xi(F)$ is generated by the
functions corresponding to loops.  The Poisson bracket is given by an
intersection pairing on oriented loops, and extended to all of
$\Xi(F)$. In \cite{quant} this is reformulated as a state sum using
unoriented loops, proving

\begin{theorem}[B-F-K]  
$K(F)$ and the map $\Phi : K_0(F \times I) \rightarrow \Xi(F)$ form a
Poisson quantization of the standard Poisson algebra $\Xi(F)$.
\end{theorem}

\section{Lattice Gauge Field Theory}

Lattice gauge field theory gives an alternative quantization of
$\Xi(F)$.  To see this, we first sketch how an $SU(2)$ gauge theory on
$F$ recovers the $SU(2)$-characters of $\pi_1(F)$.  We then pass to a
lattice model of the theory, in which a Lie group may be replaced
with its universal enveloping algebra.  Finally, the enveloping
algebra may be deformed to a quantum group.  Along the way, of course,
we will complexify to return to the $\g$ setting.

An $SU(2)$ gauge theory over $F$ consists of connections, gauge
transformations (also called the gauge group) and gauge fields.
These objects have technical definitions involving the geometry of an
$SU(2)$-bundle over $F$, but for our purposes only a few consequences
are relevant. 

First of all, a connection determines a notion of parallel transport
along a path, $\gamma$, which assigns to it an element $hol(\gamma)$
of $SU(2)$.  This element is called the holonomy of the connection
along $\gamma$. Notice that if you traverse the path in the opposite
direction then the holonomy is the inverse.  A connection is flat if
holonomy only depends on the homotopy class of a path relative to its
endpoints.

Second, the gauge group acts on connections.  A gauge
transformation can be thought of as an element of $SU(2)$ assigned to
each point of $F$.  Its effect on a connection is irrelevant; its
effect on holonomy is $hol(\gamma) \mapsto g\;hol(\gamma)\;h^{-1}$,
where $g$ and $h$ correspond to the beginning and
end points of $\gamma$.  

Finally, gauge fields are (real analytic) functions on connections.
There is an adjoint action of the gauge group on gauge
fields;\footnote{For a gauge transformation $g$, a gauge field $f$,
and a connection $x$, $(g\bullet f)(x) = f(g\bullet x)$.}
invariant gauge fields are called observables.  

Flat connections give rise to representations of $\pi_1(F)$ into $\g$
via holonomy of loops.  There are actually more flat connections than
representations. However, two connections are gauge equivalent if and
only if their holonomy representations are conjugate.

The observables, restricted to flat connections, are a space of (real
analytic) functions on $SU(2)$ representations, which are invariant
under conjugation.  The ``polynomials'' in this space---a dense
set---are the $SU(2)$-characters of $F$.

Much of the technical detail glossed over in the last few paragraphs
vanishes if we pass to a lattice model; a combinatorial setting in
which geometry is disposed of and the behavior of holonomy is
axiomatized.  As a bonus, one need not base the theory on a compact Lie
group.  What follows works for any affine algebraic group, but we will
stick to $\g$ for continuity.

Suppose that $F$ is triangulated.  The 1-skeleton of the triangulation
of $F$ is a graph.  Let $V$ denote the set of vertices and $E$ the set
of edges, each with an orientation.  The objects of a lattice gauge
field theory over $F$ are:
\begin{enumerate}
\item the connections, $\displaystyle{\bA=\prod_{e \in E} \g}$,
\item the gauge group, $\displaystyle{\cG=\prod_{v\in V} \g}$, and
\item the gauge fields, $\displaystyle{C[\bA]= \bigotimes_{e\in E}
C[\g]}$.
\end{enumerate}

\begin{figure}
\centering
\makebox[117pt]{\lmk $x_1$ \hfill $x_2$ \hfill $x_3$ \rmk}\makebox[160pt]{}\\
\raisebox{-3pt}{\mbox{\epsfig{file=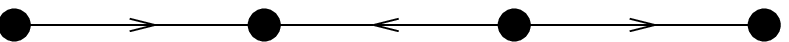,width=117pt}}}
  \makebox[160pt]{\lmk $\Longrightarrow$ \hfill 
  $hol(x_1,x_2,x_3)=x_1x_2^{-1}x_3$ \rmk}\\
\makebox[117pt]{\lmk $e_1$ \hfill $e_2$ \hfill $e_3$ \rmk}
  \makebox[160pt]{}
\caption{Example of holonomy in a lattice.}
\label{holonomy}
\end{figure}

In the formula above, $C[\g]$ is the coordinate ring of $\g$.
One thinks of a connection as assigning an element of $\g$ to each
edge.  A path is a string of edges.  Holonomy of $(x_1,x_2,x_3)$ along
the path $\{e_1,e_2,e_3\}$ is depicted in Figure \ref{holonomy}.  Note
that holonomy is clearly inverted if the path is reversed.

\begin{figure}
\centering
\makebox[30pt]{}\makebox[.15in][r]{\raisebox{-14pt}{$y$}}
                \makebox[.35in][l]{$x_2$}\makebox[90pt]{}
                \makebox[.5in]{$yx_2$}
                \makebox[30pt]{}\\ 
\vspace{-12pt}
\makebox[30pt][r]{\raisebox{.25in}{$x_3$}}
                \epsfig{file=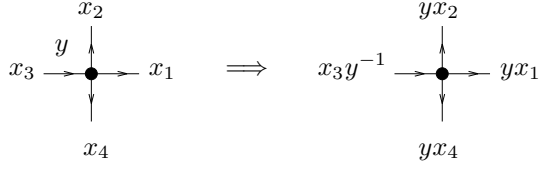,width=.5in}
                \makebox[90pt]{\raisebox{.25in}{$x_1$}\hfill
                               \raisebox{.25in}{$\Longrightarrow$}\hfill
                               \raisebox{.25in}{$x_3y^{-1}$}}
                \epsfig{file=vertex.eps,width=.5in} 
                \makebox[30pt][l]{\raisebox{.25in}{$yx_1$}}\\ 
\makebox[30pt]{}\makebox[.5in]{$x_4$}\makebox[90pt]{}
                \makebox[.5in]{$yx_4$}\makebox[30pt]{}\\
\caption{Gauge group action at a vertex.}
\label{action}
\end{figure}

One thinks of a gauge transformation as an element of $\g$ at each
vertex.  The action of the gauge group on a connection is illustrated
near a vertex in Figure \ref{action}. Note that the action is by
$y^{-1}$ on the right if an edge points in and by $y$ on the left if
it points out, a convention we adhere to through this and the next two
sections.

The gauge fields can be evaluated on connections in the obvious way.
By taking adjoints we get an action of the gauge group on the gauge
fields. The fixed subring of this action is the ring of
$\g$-characters of the one skeleton.  If $G$ is the fundamental group
of the 1-skeleton this ring is isomorphic to $\Xi(G)$.

Flatness should amount to holonomy being independent of path, but in a
lattice model we prefer the following equivalent definition.  A
connection is flat on a face of the triangulation if it is gauge
equivalent to one which has 1 on each edge of the face.  A flat
connection is flat on each face.  Invariant gauge fields evaluated on
flat connections form a ring of observables which, regardless of the
choice of triangulation, is isomorphic to
$\Xi(F)$.\footnote{Technically, divide the gauge field algebra by the
annihilator of all flat connections and then restrict to the gauge
invariant part of the quotient.}

This is an easily manipulated model of a gauge theory, but groups do
not quantize; algebras do.  So, replace $\g$ with the universal
enveloping algebra $\U$.  This is a cocommutative Hopf algebra.  The
interested reader may find a full explanation in \cite{abe} for
example, but we can get by with less.  There is an involution
$S:\U\rightarrow\U$ that corresponds to inversion in the group, a
counit $\epsilon :\U\rightarrow\bC$, and a comultiplication
$\Delta:\U\rightarrow\U\otimes\U$.  One may regard $\Delta^n$ as an
operation that breaks an element of $\U$ into states residing in
$\U^{\otimes(n+1)}$.  The notation for this is due to Sweedler
\cite{sweedler}.  For example, 
\[ \Delta^3(y) = \sum_{(y)} y^{(1)}\otimes y^{(2)}\otimes
y^{(3)}\otimes y^{(4)}. \]

\begin{figure}[b]
\centering
\makebox[30pt]{}\makebox[.5in]{\raisebox{2pt}{$e_2$}}\makebox[30pt]{}\\
\makebox[30pt][r]{\raisebox{.25in}{$e_3$}}
  \epsfig{file=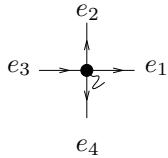,width=.5in}
  \makebox[30pt][l]{\raisebox{.25in}{$e_1$}} \\
\makebox[30pt]{}\makebox[.5in]{$e_4$}\makebox[30pt]{}
\caption{Ciliated vertex with edges ordered $e_1<e_2<e_3<e_4$.}
\label{ciliated-vertex}
\end{figure}

Since $C[\g]$ lies in the dual of $\U$, we can almost repeat the
entire process with
\begin{enumerate}
\item connections $\displaystyle{\bA=\bigotimes_{e \in E} \U}$,
\item gauge algebra $\displaystyle{\cG=\bigotimes_{v\in V} \U}$, and
\item gauge fields $\displaystyle{C[\bA]= \bigotimes_{e\in E}C[\g]}$.
\end{enumerate}
The catch is the gauge action.  In order to make sense of it, we need
to assign an ordering to the edges at each vertex.  This is done by
marking the vertex with a cilium (see Figure \ref{ciliated-vertex})
after which the orientation on $F$ gives a counter-clockwise ordering
of the edges. 

\begin{figure}
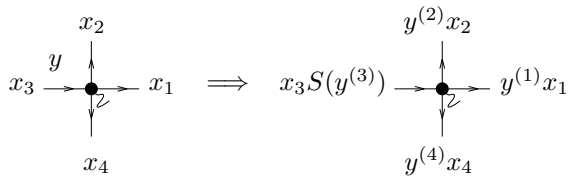

\centering
\makebox[30pt]{}\makebox[.15in]{\raisebox{-14pt}{$y$}}
                \makebox[.34in][l]{$x_2$}\makebox[90pt]{}
                \makebox[.5in]{$y^{(2)}x_2$}
                \makebox[30pt]{}\\ 
\vspace{-12pt}
\makebox[30pt][r]{\raisebox{.25in}{$x_3$}}
                \epsfig{file=ciliated-vertex.eps,width=.5in}
                \makebox[90pt]{\raisebox{.25in}{$x_1$}\hfill
                               \raisebox{.25in}{$\Longrightarrow$}\hfill
                               \raisebox{.25in}{$x_3S(y^{(3)})$}}
                \epsfig{file=ciliated-vertex.eps,width=.5in} 
                \makebox[30pt][l]{\raisebox{.25in}{$y^{(1)}x_1$}}\\ 
\makebox[30pt]{}\makebox[.5in]{$x_4$}\makebox[90pt]{}
                \makebox[.5in]{$y^{(4)}x_4$}\makebox[30pt]{}\\
\caption{Gauge algebra action at a vertex.}
\label{q-action}
\end{figure}

It is best to think of connections and gauge transformations as pure
tensors, remembering always to extend linearly.  We thus view a
connection as an assignment of an element of $\U$ to each edge of the
triangulation.  Holonomy is apparent; for the path in Figure
\ref{holonomy} it would be $x_1S(x_2)x_3$, where $S$ is the antipode
of $\U$.  We continue to think of a gauge transformation as an element
of $\U$ at each vertex, with the action at a vertex illustrated in
Figure \ref{q-action}.

A further problem with the action is that gauge ``equivalence'' is not
an equivalence relation anymore, necessitating a slight technical
modification of flatness which we will not address here.  Also, the
word ``invariant'' means $y\bullet x = \epsilon(y) x.$ However, the
passage to gauge fields on flat connections modulo the gauge algebra
proceeds as before, giving exactly the same ring.

Finally we pass to $\Uh$.  This is a quasi-triangular ribbon Hopf
algebra \cite{kassel}.  It is non-cocommutative in a fashion
constrained by an element of $\Uh\otimes\Uh$ called the universal
$R$-matrix.  The antipode $S$ is no longer an involution; rather $S^2$
acts as conjugation by the so-called charmed element, $k$.  The
definition of flat connection is further altered, preserving
independence of path but deforming the holonomy of a trivial loop to
$k^{\pm 1}$.  The dual of $\Uh$ contains a deformation, $\B$, of
$C[\g]$.  Thus one hopes to obtain a quantized ring of observables by
replacing each object with its quantum analogue.

There is one small problem.  The natural multiplication on
$C[\bA]= \otimes_{e\in E}\B$ (i.e. the one dual to
the natural comultiplication on $\bA=\otimes_{e \in
E} \Uh$) is not gauge invariant.  This is a major obstruction, and
the solution is notable enough to occupy the next section.  However,
once it has been addressed, we will have

\begin{theorem} 
Quantum observables exist.  They form a ring,
$\Xi_h(F)$, which is independent of triangulation and ciliation, and
which quantizes $\Xi(F)$.  
\end{theorem}

\section{Nabla}

\begin{figure}[b]
\centering
\makebox[30pt]{}\makebox[.5in]{\raisebox{2pt}{$e_2$}}\makebox[30pt]{}\\
\makebox[30pt][r]{$e_3$}
  \epsfig{file=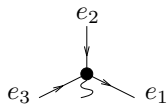,width=.5in}
  \makebox[30pt][l]{$e_1$}
\caption{Trivalent ciliated vertex.}
\label{trivalent}
\end{figure}

The natural comultiplication on the coalgebra of quantum connections
is a tensor power of $\Delta$ composed with a permutation.  For
instance,  it would send 
\[x_1\otimes x_2 \mapsto x_1^{(1)}\otimes x_2^{(1)}\otimes 
                                x_1^{(2)}\otimes x_2^{(2)}.\]
Expanding on a theme of quantum topology, we denote this morphism by a
tangle built from branches for each application of $\Delta$ and a
braid corresponding to the permutation.  We then obtain a quantized
comultiplication
\[\nabla : \bA \rightarrow \bA\otimes \bA \]
by allowing crossings to encode actions of the $R$-matrix.

\begin{figure}
\centering
\epsfig{file=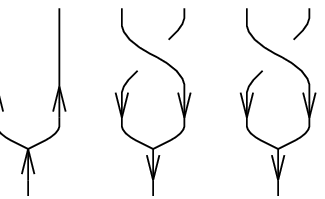,width=1.5in,height=.5in}\\
$e_1$\hspace{.5in}$e_2$\hspace{.5in}$e_3$
\caption{Coupons for each edge.}
\label{level-one}
\end{figure}

\begin{figure}
\centering
\epsfig{file=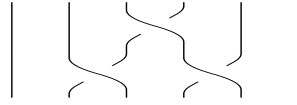,width=1.5in,height=.5in}
\caption{Six-braid encoding the permutation $(1)(2453)(6)$.}
\label{level-two}
\end{figure}

There is a fundamental tangle associated to any vertex---the one in
Figure \ref{trivalent}, for example---whose construction proceeds in
stages.  First, assign a coupon to each edge as in Figure
\ref{level-one}.
There are two types of these, depending on whether the edge points in
or out, and they must be ordered left to right matching the cilial
order of the edges.  Next, we construct a $2n$-braid ($n=\text{valence
of the vertex}$) by dragging odd numbered strands left and even
numbered strands right.\footnote{The inherent ambiguity evaporates
when we construct the morphism because the $R$-matrix solves the
Yang-Baxter equation \cite{kassel}.  It is an elegant feature of
quantum topology that isotopies of tangles correspond to identities in
a quantum group.} Evens lie over odds.  Our example is the the 6-braid
in Figure \ref{level-two}.
The fundamental tangle is formed by stacking the braid atop the
coupons.  Orientation of the coupons carries over to the strands of
the braid.

\begin{figure}[b]
\centering
\makebox[.4in]{$x''$}\makebox[.4in]{$x'$}\hspace{.4in}
  \makebox[.4in]{$x'$}\makebox[.4in]{$x''$}\\
\vspace{2pt}
\epsfig{file=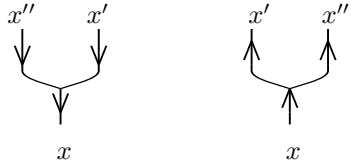,width=1.6in,height=.5in}\\
\makebox[.8in]{$x$}\hspace{.4in}\makebox[.8in]{$x$}
\caption{Comultiplication acting at a branch.}
\label{branches}
\end{figure}

Now imagine $x_1\otimes x_2\otimes x_3$ entering the tangle from the
bottom and traveling upward.  Each branch indicates comultiplication
with the output ordered as in Figure \ref{branches}.
Note that we are suppressing the summation symbols.  

\begin{figure}
\centering
\makebox[.4in]{$\beta_iy$}\makebox[.4in]{$\alpha_ix$}\hspace{.4in}
 \makebox[.4in]{$\beta_iy$}\makebox[.4in]{$xS(\alpha_i)$}\hspace{.4in}
 \makebox[.4in]{$yS(\beta_i)$}\makebox[.4in]{$\alpha_ix$}\hspace{.2in}
 \makebox[.5in]{$yS(\beta_i)$}\makebox[.5in]{$xS(\alpha_i)$} \\
\vspace{2pt}
\epsfig{file=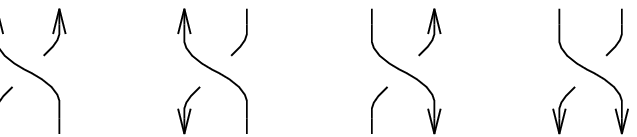,width=4in,height=.5in}\\
\makebox[.4in]{$x$}\makebox[.4in]{$y$}\hspace{.4in}
     \makebox[.4in]{$x$}\makebox[.4in]{$y$}\hspace{.4in}
     \makebox[.4in]{$x$}\makebox[.4in]{$y$}\hspace{.2in}
     \makebox[.5in]{$x$}\makebox[.5in]{$y$}
\caption{Actions of the $R$-matrix at crossings.}
\label{crossings}
\end{figure}

Each crossing corresponds to an action of the $R$-matrix, which we
write as $R = \sum_i \alpha_i\otimes \beta_i$.  The four possibilities
are shown in Figure \ref{crossings}, again suppressing summation.
Note that, as usual, left and right multiplication correspond
respectively to outward and inward pointing edges, and that right
multiplication is preceded by an application of $S$.

Sweeping out, we get a morphism
\[ \bF_v : \bigotimes_{\text{edges at $v$}} \Uh \rightarrow
   \left(\bigotimes_{\text{edges at $v$}} \Uh\right)^{\otimes 2}.\]
In our example,
\begin{multline*}
x_1\otimes x_2 \otimes x_3 \mapsto \\ 
 x_1'\otimes x_2'S(\beta_1)S(\beta_3) \otimes
 x_3'S(\beta_2)S(\beta_4)S(\beta_5) \otimes \alpha_5\alpha_3 x_1''
 \otimes x_2''S(\alpha_2)S(\alpha_4) \otimes x_3''S(\alpha_2).
\end{multline*} 
Eight summations are suppressed, and the subscripts on $\alpha_j$ and
$\beta_j$ are shorthand for summation over the $j$-th application of
the $R$-matrix.

The morphism $\bF_v$ is coassociative in the sense that $(Id \otimes
\bF_v)\circ \bF_v = ( \bF_v\otimes Id)\circ \bF_v$.  Furthermore, its
effect in a given factor of $(\bigotimes_{\text{edges at $v$}}
\Uh)^{\otimes 2}$ is either entirely by right multiplication or
entirely by left multiplication.  This allows us to combine the
effects of $\{\bF_v\;|\;v\in V\}$ into a single morphism
\[\nabla : \bA \rightarrow \bA \otimes \bA.\]

\begin{theorem}
$\nabla$ is coassociative and gauge invariant.
\end{theorem}

\section{Quantum Observables for $\g$}

\begin{figure}[b]
\centering
\makebox[2.5in]{
\lmk $v_4$ \hfill \raisebox{-.2in}{$e_4$} \hfill
\raisebox{-.4in}{$v_3$} \hfill \raisebox{-.2in}{$e_1$} \hfill $v_1$
\rmk} \\ \vspace{-.3in} 
\makebox[2.5in]{
\raisebox{.5in}{$e_5$} \hfill \epsfig{file=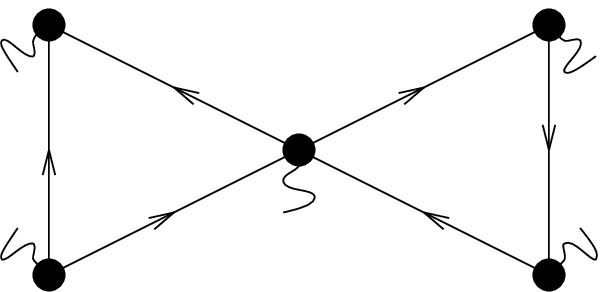,width=2in} 
\hfill \raisebox{.5in}{$e_2$}} \\ \vspace{-.1in}
\makebox[2.5in]{
\lmk $v_5$ \hfill \raisebox{.2in}{$e_6$} \hfill \mbox{} 
\hfill \raisebox{.2in}{$e_3$} \hfill $v_2$ \rmk}
\caption{Oriented ciliated graph $\Gamma$.}
\label{bowtie}
\end{figure}

At the end of Section 2 we saw how loops became functions generating
$\Xi(F)$.  In this section we will describe the quantum analogue of
that fact.  All definitions are given in terms of a running example,
so assume that $\Gamma$ is the oriented, ciliated graph in Figure
\ref{bowtie}.

Following Section 6 we have

\begin{enumerate}
\item the connection coalgebra, 
$\displaystyle{\bA =\Uh^{\otimes 6}}$,
\item the gauge algebra, 
$\displaystyle{\cG=\Uh^{\otimes 5}}$, and
\item the algebra of gauge fields, 
$\displaystyle{C[\bA] = (\B)^{\otimes 6}}$.
\end{enumerate}

Bowing to technicalities, a loop in $\Gamma$ will be allowed to meet
each edge at most once, and each vertex at most twice. In accordance
with the theme of quantization by crossings, we say a $q$-loop is a
loop with a choice of under or over crossing whenever it intersects
itself transversely.  We express this as a sequence of edges with $+$
and $-$ signs interspersed.  For example, 
\[l = \{e_1,e_2,e_3,+,e_4,e_5,e_6,-\} \]
is a $q$-loop.  It defines an element of $C[\bA]$ via the following
graphical recipe.
 
\begin{enumerate}

\item Choose a pure tensor $x = x_1\otimes\cdots\otimes x_6 \in \bA$.
Draw a picture of $\Gamma$ with $x_i$ labeling each corresponding edge.
\item Apply $\epsilon$ to any edge not appearing in the loop.  (No effect
  on this example.)  
\item At the crossing, act by an $R$-matrix, either
$\sum\alpha_i\otimes\beta_i$ or $\sum\beta_i\otimes S(\alpha_i)$.  The
action will take place on the first two edges in the ciliation; the
$R$-matrix is chosen so that $\beta_i$ acts on the bottom strand
($\sum\alpha_i\otimes\beta_i$ in this example); and the action follows
left/right rules as in Section 5.  Write this on the appropriate
edges, suppressing summations.
\item If an edge is oriented against the direction of the loop,
multiply on the right by $k$, and then apply $S$.  
\item Each time the loop passes through a vertex check to see if the
incoming edge goes before the outgoing edge in the ciliation.  If not,
then right multiply by $k$ on the incoming edge.  The picture should
now look like Figure \ref{sample-loop}.
\begin{figure}
\centering
\makebox[1in]{}
  \makebox[1in][l]{\hspace{.5in}$x_4k$}
  \makebox[1in][r]{$\beta_ix_1$\hspace{.5in}}
  \makebox[1in]{} \\ \vspace{-.2in}
\makebox[1in][r]{\raisebox{.5in}{$S(x_5k)k$}}\makebox[2in]{
  \epsfig{file=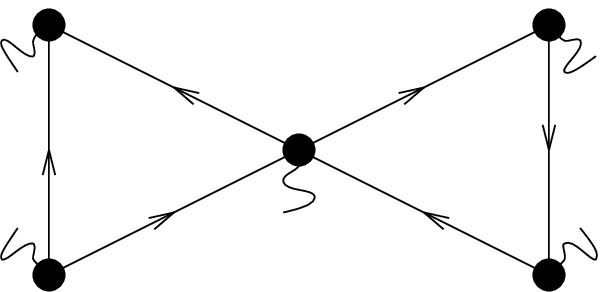,width=2in}} 
  \makebox[1in][l]{\raisebox{.5in}{$x_2$}} \\ \vspace{-.2in}
\makebox[1in]{}
  \makebox[1in][l]{\hspace{.5in}$x_6k$}
  \makebox[1in][r]{$x_3S(\alpha_i)$\hspace{.5in}}
  \makebox[1in]{} 
\caption{Action of the $q$-loop $l$ on $x$.}
\label{sample-loop}
\end{figure}
\item Multiply everything together as you traverse the loop. Take the
image of this ``quantum holonomy'' in the fundamental representation
(see \cite{KM} or \cite{RT1}) of $U_h(sl_2)$. Finally, take the trace
to get a complex number:
\[x \mapsto \sum_i\tr(\beta_ix_1x_2x_3S(\alpha_i)x_4S(x_5)kx_6k).\]
\item Extend linearly over all of $\bA$ to obtain a function
$W_l$. 
\end{enumerate}

The function we have defined is usually called a Wilson loop in the
literature.  It should be clear that a $q$-loop is just a knot diagram
with a base point and an orientation.  Our goal is to assign a quantum
observable to each equivalence class of link diagrams.  

The first step is to note that the rules we gave for acting on edges
prior to computing holonomy are local.  One could just as well apply
them to a link of loops, compute holonomy along each, and take the
product of the resulting traces.  Clearly the individual traces are
independent of base points.  Reversing orientation is less trivial, as
it involves $S$ rather than inversion.  But $\tr(S(x))=\tr(x)$ in the
fundamental representation, so orientations don't matter.  Gauge
invariance is easily checked at individual vertices. Finally, suppose
that $\Gamma$ is the 1-skeleton of a triangulated surface. Let $l$ and
$l'$ be equivalent link diagrams and $x$ a flat connection.  Since
flatness implies independence of path, $W_l = W_{l'}$.

At this point we see that a link $L$ determines an observable, $W_L$,
provided one has a fine enough triangulation of $F$.  Since $\Xi_h(F)$
is independent of triangulation, we can make the assignment
\[ L \mapsto (-1)^{|L|}W_L,\]
where $|L|$ is the number of components of $L$.  Linearity and
continuity extend it uniquely to a map
\[\Phi_h : \bC\cL_{F\times I}[[h]] \rightarrow \Xi_h(F).\]  

This is the quantum analogue of $\Phi$ from Section 2, taking loops to
the character ring.  That map took (adjoints of) the skein relations
in $K_0$ to the fundamental $\g$-trace identity and to $\tr(I)=2$.
The corresponding quantized identities in $\Uh$ are
\begin{align*}
t\;\tr(ZY)+t^{-1}\;\tr(S(Z)W) &= \sum_i\tr(\alpha_iz)\tr(\beta_ix),
   \quad\text{and}\\
\tr(k^{\pm 1})& = t^2+t^{-2}. 
\end{align*}

It is clear from the definition of flat connection that the (adjoint
of the) skein relation $\bigcirc+(t^2+t^{-2})$ maps to $\tr(k) =
t^2+t^{-2}$, but
\[ \lcr+t\zer+t^{-1}\ift\]
is more complicated because it has less symmetry than the quantum
trace identity.  In some cases, the (adjoint) skein relation is
obviously mapped to an identity, while others require more manipulation.

\begin{theorem}
$\Phi_h(S(F)) = 0$. Furthermore, the quotient map 
\[\widetilde{\Phi}_h:K(F) \rightarrow \Xi_h(F)\]
is an isomorphism.
\end{theorem}

\section{The Future}

The results described in this survey place skein theory at the
confluence of ideas from topology, representation theory,
noncommutative algebra and mathematical physics. Standard
techniques from skein theory \cite{HP3,HP4} extend our lattice
construction of  $K(F)$ to a description of
the Kauffman bracket skein module of an arbitrary compact
$3$-manifold.  Consequently, $K(M)$ has an intensive definition in
terms of links and skein relations, and an extensive definition in
terms of quantized invariant theory.  This nexus suggests some avenues
for further research.

It is proved in \cite{quant} that the affine $\g$-characters induce
topological generators of $K(M)$. In particular, the Kauffman bracket
skein module of a small $3$-manifold (i.e. containing no
incompressible surface) is finitely generated, and thus can be used as
a classification tool.

If $K(M)$ is topologically free then there is a meaningful 
pairing between it and the set of equivalence classes of
$\g$-representations of $\pi_1(M)$. For nilpotent free $K_0(M)$ this
is a duality pairing. In the case of $F\times I$ the pairing has an
especially easy form because the basis is a canonical set of links
\cite{quant}.  The Yang-Mills measure on the algebra of observables
can be computed along the same lines \cite{Bu}. This holds out the
promise of producing integral formulas for the
Witten-Reshetikhin-Turaev invariants of a $3$-manifold that will admit
to asymptotic analysis.

The focus in this paper has been $\g$, but the lattice gauge field
theory works for any algebraic group \cite{lattice}.  There should be
skein modules corresponding to the other groups, just as $K(M)$
corresponds to $\g$. We will need two kinds of skein relations:
fundamental relations in the Hecke algebra associated to the group,
and the quantized Cayley-Hamilton identity. There has been some study
of these ideas due to Kuperberg \cite{spiders} and Anderson, Mattes
and Reshetikhin \cite{AMR}.

In another direction, it should be possible to commence the study of
the syzygies of skein modules. A syzygy is a relationship between
relationships. For instance, we can define a homology theory for the
Kauffman bracket skein module.  The 0-chains are spanned by all links;
the 1-chains by all ``Kauffman bracket skein triples''; the 2-chains
by all ``triples of triples'', etc.  The 0-th homology of this complex
is $K(M)$, and have examples to show that the theory is not always
trivial.  Notice that the $n$-th homology is measuring relations among
relations.

The opacity of the structure of $K(M)$ poses many questions. Is it
possible for $K(M)$ of a compact manifold to have torsion and still be
topologically finitely generated?  What is the relationship between
$K_A(M)$ and $K(M)$? Przytycki has an example of a noncompact manifold
where $K_A(M)$ is infinitely generated yet $K(M)$ is trivial. There is
a grading of $K(M)$ by cables.  The top term in the grading is
everything; after that you take the span of all $2$-fold cables, then
$3$-fold cables, etc.  How is torsion in $K(M)$ reflected in this grading?
Are there nilpotents in any $K_0(M)$, and if so, how do they affect
the geometry of the representation space?  Finally, in the interest of
computability, what is a relative skein module and is there a gluing
theorem?

The Kauffman bracket skein module is organic to many fields.  We hope
it, and other skein modules, will act as  catalysts for the
synergistic mixing of ideas from these fields.

\end{document}